\documentclass[11pt,a4paper]{amsart} 
\usepackage[all]{xy}
\SelectTips{cm}{}
\usepackage{ifthen} 
\usepackage{amssymb}
\usepackage{stmaryrd}
\calclayout
\makeatletter 
\def\serieslogo@{} 
\def\@setcopyright{} 
\makeatother

\title{Thick subcategories and virtually Gorenstein algebras}

\author{Apostolos Beligiannis}
\address{Apostolos Beligiannis\\Department of Mathematics\\ 
University of Ioannina\\45110 Ioannina\\Greece.}
\email{abeligia@cc.uoi.gr}
\author{Henning Krause}
\address{Henning Krause\\ Institut f\"ur Mathematik\\
Universit\"at Paderborn\\ 33095 Paderborn\\ Germany.}
\email{hkrause@math.uni-paderborn.de}
\thanks{Version from August 29, 2006.}

\newtheorem*{lem}{Lemma}
\newtheorem*{prop}{Proposition} 
\newtheorem*{cor}{Corollary}
\newtheorem*{thm}{Theorem} 

\newtheorem*{Thm1}{Theorem~1}
\newtheorem*{Thm2}{Theorem~2}

\theoremstyle{remark}

\theoremstyle{definition}
\newtheorem*{exm}{Example}
\newtheorem*{defn}{Definition}
\newtheorem*{rem}{Remark}
\numberwithin{equation}{section}

\renewcommand{\mod}{\operatorname{mod}\nolimits}

\newcommand{\Thick}{\operatorname{Thick}\nolimits}
\newcommand{\Thickc}{\Thick^{\scriptscriptstyle\coprod}}
\newcommand{\Thickp}{\Thick^{\scriptscriptstyle\prod}}
\newcommand{\proj}{\operatorname{proj}\nolimits}
\newcommand{\Proj}{\operatorname{Proj}\nolimits}
\newcommand{\Inj}{\operatorname{Inj}\nolimits}
\newcommand{\GProj}{\operatorname{GProj}\nolimits}
\newcommand{\Gproj}{\operatorname{Gproj}\nolimits}
\newcommand{\GInj}{\operatorname{GInj}\nolimits}
\newcommand{\Ginj}{\operatorname{Ginj}\nolimits}
\newcommand{\uGProj}{\operatorname{\underline{GProj}}\nolimits}
\newcommand{\oGInj}{\operatorname{\overline{GInj}}\nolimits}
\newcommand{\inj}{\operatorname{inj}\nolimits}
\newcommand{\rad}{\operatorname{rad}\nolimits}

\newcommand{\Mod}{\operatorname{Mod}\nolimits}

\newcommand{\Hom}{\operatorname{Hom}\nolimits}

\newcommand{\Ker}{\operatorname{Ker}\nolimits}
\newcommand{\Coker}{\operatorname{Coker}\nolimits}

\newcommand{\Ext}{\operatorname{Ext}\nolimits}

\newcommand{\umod}{\operatorname{\underline{mod}}\nolimits}
\newcommand{\uMod}{\operatorname{\underline{Mod}}\nolimits}
\newcommand{\oMod}{\operatorname{\overline{Mod}}\nolimits}

\newcommand{\oHom}{\operatorname{\overline{Hom}}\nolimits}

\newcommand{\omod}{\operatorname{\overline{mod}}\nolimits}

\newcommand{\tac}{\mathrm{tac}}
\newcommand{\op}{\mathrm{op}}

\newcommand{\lto}{\longrightarrow}
\newcommand{\xto}{\xrightarrow}

\def\li{\varinjlim}

\def\a{\alpha}

\def\d{\delta}

\def\p{\phi}

\def\La{\Lambda}

\def\A{{\mathcal A}}

\def\C{{\mathcal C}}
\def\D{{\mathcal D}}
\def\E{{\mathcal E}}
\def\I{{\mathcal I}}

\def\X{{\mathcal X}}
\def\Y{{\mathcal Y}}

\def\T{{\mathcal T}}
\def\U{{\mathcal U}}
\def\V{{\mathcal V}}

\def\bbZ{\mathbb Z}

\def\bfi{\mathbf i}
\def\bfp{\mathbf p}
\def\bfs{\mathbf s}

\def\bfD{\mathbf D}
\def\bfK{\mathbf K}

\begin{document}

\begin{abstract}
An Artin algebra is by definition virtually Gorenstein if the class of
modules which are right orthogonal (with respect to $\Ext^*(-,-)$) to
all Gorenstein projective modules coincides with the class of modules
which are left orthogonal to all Gorenstein injective modules. We
provide a new characterization in terms of finitely generated
modules. In addition, an example of an algebra is presented which is
not virtually Gorenstein.
\end{abstract}
\maketitle

\section{Introduction}

Let $\La$ be an Artin algebra.  We consider the category $\Mod\La$ of
right $\La$-modules and its full subcategory $\mod\La$ of finitely
generated modules. Let $\Proj\La$ denote the full subcategory of
$\Mod\La$ consisting of all projective modules and let
$\proj\La=\Proj\La\cap\mod\La$.  Analogously the subcategories of
injective $\La$-modules $\Inj\La$ and $\inj\La$ are defined.

Given any exact category $\A$, we call a full subcategory $\C$ of
$\A$ {\em thick} if it is closed under direct factors and has the
following two out of three property: for every exact sequence $0\to
X\to Y\to Z\to 0$ in $\A$ with two terms in $\C$, the third term
belongs to $\C$ as well. Let $\Thick(\C)$ denote the smallest thick
subcategory of $\A$ which contains $\C$. Take for example
$\A=\mod\La$. Then $\Thick(\proj\La)$ consists of all finitely
generated $\La$-modules having finite projective dimension and
$\Thick(\inj\La)$ consists of all finitely generated $\La$-modules
having finite injective dimension.

The algebra $\La$ is by definition {\em right Gorenstein} if $\La$
viewed as a right module has finite injective dimension, equivalently
if $\Thick(\proj\La)\subseteq \Thick(\inj\La)$. Note that it is an
open problem, whether any left Gorenstein algebra is right
Gorenstein. More specifically, in \cite{AR} 
the following connection with the finitistic dimension conjecture is pointed out.
\medskip

{\em For a left Gorenstein algebra
$\La$ the following conditions are equivalent. 
\begin{enumerate}
\item The algebra $\La$ is right Gorenstein.
\item The finitistic dimension of $\La$ is finite.
\item The subcategory $\Thick(\proj\La)$ of $\mod\La$ is contravariantly finite.
\end{enumerate}}
\medskip

Recall that a subcategory $\X$ of $\mod\La$ is {\em contravariantly
finite} if every finitely generated $\La$-module $C$ has a {\em right
$\X$-approximation} $X\to C$, that is, $X\in\X$ and the induced map
$\Hom_\La(X',X)\to\Hom_\La(X',C)$ is surjective for every
$X'\in\X$. {\em Covariantly finite} subcategories are defined dually.

The problem to understand the Gorenstein left-right symmetry provides
a motivation for studying the following class of algebras which has
been introduced in \cite{BR}.  An algebra $\La$ is {\em virtually
Gorenstein} if for every $\La$-module $X$, the functor
$\Ext^i_\La(X,-)$ vanishes for all $i>0$ on all Gorenstein injective
$\La$-modules if and only if $\Ext^i_\La(-,X)$ vanishes for all $i>0$
on all Gorenstein projective $\La$-modules; see (\ref{se:vgor}) for
details. In this note we provide the following characterization in
terms of finitely generated modules, complementing the discussion of
virtually Gorenstein algebras in \cite{B}.

\begin{Thm1}
For an Artin algebra $\La$ the following are equivalent.
\begin{enumerate}
\item The algebra $\La$ is virtually Gorenstein.
\item The subcategory $\Thick(\proj\La\cup\inj\La)$ of $\mod\La$ is
contravariantly finite.
\item The subcategory $\Thick(\proj\La\cup\inj\La)$ of $\mod\La$ is
covariantly finite.
\end{enumerate}
\end{Thm1} 

An immediate consequence of this result is the fact that $\La$ is
virtually Gorenstein if and only if $\La^\op$ is virtually
Gorenstein. The proof of the theorem relies on the interplay between
cotorsion pairs for modules and torsion pairs for complexes of
modules, building on previous work of both authors \cite{B,K1}.  The
crucial step is a description of the subcategory
$\Thick(\proj\La\cup\inj\La)$ which we formulate as a separate result.

\begin{Thm2}
For a finitely generated $\La$-module $X$ the following are equivalent.
\begin{enumerate}
\item The module $X$ belongs to $\Thick(\proj\La\cup\inj\La)$.
\item $\Ext^i_\La(X,-)$ vanishes for all $i>0$ on all Gorenstein injective
$\La$-modules.
\item $\Ext^i_\La(-,X)$ vanishes for all $i>0$ on all Gorenstein
projective $\La$-modules.
\end{enumerate}
\end{Thm2}

\section{Gorenstein projective and injective modules}
We recall the definitions of Gorenstein projective and Gorenstein
injective modules. These classes of modules induce two cotorsion pairs
which give rise to the property of an algebra to be virtually
Gorenstein.

\subsection{Cotorsion pairs for abelian categories}
Let $\A$ be an abelian category, for instance $\A=\Mod\La$ or
$\A=\mod\La$.  Let $\X$ and $\Y$ be classes of objects in $\A$. We
define
$$\X^\perp =\{Y\in\A\mid
\mbox{$\Ext^i_\A(X,Y)=0$ for all $i>0$ and $X\in\X$}\},$$
$$^\perp \Y=\{X\in\A\mid \mbox{$\Ext^i_\A(X,Y)=0$ for all $i>0$ and
$Y\in\Y$}\}.$$ Now let $$0\lto A\stackrel{\p}\lto B
\stackrel{\psi}\lto C\lto 0$$ be an exact sequence in $\A$.  The map
$\psi$ is called a {\em special right $\X$-approximation} of $C$ if
$B\in\X$ and $A\in\X^\perp$.  Dually, the map $\p$ is called a {\em
special left $\Y$-approximation} of $A$ if $B\in\Y$ and $C\in
{^\perp\Y}$.  The following concept of a cotorsion pair is due to
Salce.

\begin{defn}\label{de:cotorsion}
A {\em cotorsion pair} for $\A$ is a pair $(\X,\Y)$ of subcategories
of $\A$ satisfying the following conditions:
\begin{enumerate}
\item $\X={^\perp\Y}$ and $\Y=\X^\perp$;
\item every object in $\A$ admits a special right $\X$-approximation
and a special left $\Y$-approximation.
\end{enumerate}
\end{defn}

\subsection{Gorenstein projective and Gorenstein injective modules}\label{se:vgor}

Let $\A$ be an additive category.  A complex
$$X=\;\;\cdots\lto X^{i-1}\lto X^{i}\lto X^{i+1}\lto\cdots$$ in
$\A$ is called {\em totally acyclic} if $\Hom_\A(A,X)$ and
$\Hom_\A(X,A)$ are acyclic complexes of abelian groups for all $A$ in
$\A$.  The following classes of modules have been introduced by
Auslander.

\begin{defn}
A $\La$-module $C$ is called
\begin{enumerate} 
\item {\em Gorenstein projective} if it is of
the form $C=\Coker (X^{-1}\to X^0)$ for some totally acyclic complex
$X$ of projective $\La$-modules, and
\item {\em Gorenstein injective} if it is of the fom $C=\Ker
(X^{0}\to X^1)$ for some totally acyclic complex $X$ of injective
$\La$-modules.
\end{enumerate}
\end{defn}

Note that different terminology is used in the literature. We have 
$$\text{Gorenstein projective = G-dimension zero =
Cohn-Macaulay.}$$ We denote by $\GProj\La$ the full subcategory of
$\Mod\La$ which is formed by all Gorenstein projective $\La$-modules,
and $\GInj\La$ denotes the full subcategory which is formed by all
Gorenstein injective $\La$-modules.

In \cite{BR}, it is shown that there are cotorsion pairs
$$(\GProj\La,(\GProj\La)^\perp) \quad\text{and} \quad
({^\perp(\GInj\La)},\GInj\La)$$ for $\Mod\La$ satisfying
$$\GProj\La\cap(\GProj\La)^\perp=\Proj\La \quad\text{and} \quad
{^\perp(\GInj\La)}\cap\GInj\La=\Inj\La.$$
The algebra $\La$ is called {\em virtually Gorenstein} if 
$$(\GProj\La)^\perp={^\perp(\GInj\La)}.$$ We refer to \cite{B} for an
extensive discussion of virtually Gorenstein algebras. For example, a
left Gorenstein algebra is virtually Gorenstein if and only if it is
right Gorenstein. In that case we have
$$(\GProj\La)^\perp=\Thick(\Proj\La)=\Thick(\Inj\La)={^\perp(\GInj\La)}.$$
Also, every algebra of finite representation type is virtually
Gorenstein, and every algebra derived equivalent to a virtually
Gorenstein algebra is again virtually Gorenstein.

\subsection{Resolving and coresolving subcategories}\label{se:resolv}
We include for later reference some basic facts about resolving
subcategories.  Recall that a subcategory $\X$ of an abelian category
$\A$ is {\em resolving} if it contains all projectives, and if for
every exact sequence $0\to A\to B\to C\to 0$ in $\A$, we have
$A,C\in\X$ implies $B\in\X$, and $B,C\in\X$ implies $A\in\X$.  {\em
Coresolving} subcategories are defined dually.

\begin{lem}[{\cite[Theorem~7.12]{B}}]
Let $(\X,\Y)$ be a cotorsion pair for $\Mod\La$. 
\begin{enumerate}
\item The subcategory $\X$ is resolving.
\item The subcategory $\Y$ is resolving if and only
$\X\cap\Y=\Proj\La$.
\item If $\Y$ is resolving, then $\X\subseteq\GProj\La$.
\end{enumerate}
\end{lem}
\begin{proof}
(1) is clear. To prove (2), suppose first that $\Y$ is
resolving. Clearly, $\Proj\La\subseteq\X\cap\Y$. Now let $Z\in\X\cap
\Y$ and choose an epimorphism $\p\colon P\to Z$ with $P$
projective. We have $\Ker\p\in\X\cap\Y$ and therefore $\p$
splits. Thus $Z$ is projective and we have $\X\cap\Y=\Proj\La$.  Now
suppose that $\X\cap\Y=\Proj\La$. Then every epimorphism $U\to V$ with
$U,V\in \X\cap\Y$ splits. This implies that $\Y$ is a thick
subcategory; see \cite[Lemma~3.2]{KS}. In particular, $\Y$ is resolving.

To prove (3), suppose again that $\Y$ is resolving and fix
$X\in\X$. We choose a projective resolution
$$\cdots \lto P^{-3}\stackrel{\d^{-3}}\lto
P^{-2}\stackrel{\d^{-2}}\lto P^{-1}\stackrel{}\lto X\lto 0$$ of
$X$ and complete this inductively to a totally acyclic complex
$(P,\d)$ of projective modules as follows. Choose for each $i\geq
-1$ a special left $\Y$-approximation $\Coker\d^{i-1}\to P^{i+1}$ and
take for $\d^i$ the composition $P^{i}\to\Coker\d^{i-1}\to P^{i+1}$.
Note that each $P^i$ is projective since a special left
$\Y$-approximation of an object in $\X$ belongs to $\X\cap\Y$.
\end{proof}

We mention that there is a dual result about coresolving subcategories
and Gorenstein injective modules.

\section{Categories of complexes and torsion pairs}

\subsection{The category $\bfK(\Inj\La)$}\label{se:kinj}

Let $\A$ be an additive category. We denote by $\bfK(\A)$ the homotopy
category of complexes in $\A$. If $\A$ is abelian, then $\bfD(\A)$
denotes the derived category of $\A$ which is obtained from $\bfK(\A)$
by formally inverting all quasi-isomorphisms. Given any category $\T$
with small coproducts, we denote by $\T^c$ the full subcategory which
is formed by all compact objects. Recall that an object $X$ is {\em
compact} if $\Hom_\T(X,-)$ preserves small coproducts, equivalently if
every map $X\to\coprod_iY_i$ into a small coproduct factors through a
finite coproduct of $Y_i$s.

\begin{prop}[{\cite[Proposition~2.3]{K1}}]\label{pr:KInj}
The triangulated category $\bfK(\Inj\La)$ is compactly generated.  The
canonical functor $\bfK(\Inj\La)\to\bfD(\Mod\La)$ induces an
equivalence
$$\bfK(\Inj\La)^c\stackrel{\sim}\lto\bfD^b(\mod\La).$$
\end{prop}

We shall use the functor
$$\bfi\colon\Mod\La\lto\bfK(\Inj\La)$$ which sends a $\La$-module to
its injective resolution. More generally, the canonical functor
$\bfK(\Inj\La)\to\bfD(\Mod\La)$ has a right adjoint sending a complex
$X$ to its semi-injective resolution $\bfi X$.  This right adjoint
induces an equivalence $\bfD^b(\mod\La)\xto{\sim}\bfK(\Inj\La)^c$
which is a quasi-inverse for the equivalence
$\bfK(\Inj\La)^c\xto{\sim}\bfD^b(\mod\La)$.

The following functor takes complexes to modules. We define
$$\bfs\colon\bfK(\Inj\La)\lto\oMod\La,\quad X\mapsto\Ker (X^0\to
X^1),$$ where $\oMod\La$ denotes the {\em stable category modulo
injectives}.  The objects of $\oMod\La$ are those of $\Mod\La$ and the
morphisms for a pair $X,Y$ of $\La$-modules are by definition
$$\oHom_\La(X,Y)=\Hom_\La(X,Y)/\I(X,Y),$$ where $\I(X,Y)$ denotes the
subgroup of morphisms which factor through an injective
module. Analogously, the full subcategories $\oGInj\La$ and $\omod\La$ of
$\oMod\La$, and the stable category $\uMod\La$ modulo projectives are defined.

\begin{rem}
  The category $\bfK(\Proj\La)$ is compactly generated, simply because
  the equivalence $\Proj\La\xto{\sim}\Inj\La$ induces an equivalence
  $\bfK(\Proj\La)\xto{\sim}\bfK(\Inj\La)$. However, in general it is not true
  that the projective resolution $\bfp X$ of a finitely generated
  $\La$-module $X$ is a compact object in $\bfK(\Proj\La)$. This is
  precisely the reason for concentrating on the category
  $\bfK(\Inj\La)$ in this work.
\end{rem}

\subsection{Totally acyclic complexes}\label{se:tac}

We denote by $\bfK_\tac(\A)$ the full subcategory of $\bfK(\A)$ which
is formed by all totally acyclic complexes.  The following proposition
is the analogue of a result for the category $\bfK(\Proj\La)$ due to
J{\o}rgensen \cite{J}.

\begin{prop}[{\cite[\S 7]{K1}}]
\begin{enumerate}
\item The functor $\bfs$ induces an equivalence
$\bfK_\tac(\Inj\La)\xto{\sim} \oGInj\La$.
\item The inclusion functor $\bfK_\tac(\Inj\La)\to\bfK(\Inj\La)$ has a
left adjoint which induces (via $\bfs$) a left adjoint $F\colon
\oMod\La\to\oGInj\La$ of the inclusion $\oGInj\La\to\oMod\La$.
\item Given $X\in\Mod\La$, the adjunction morphism $X\to FX$
represents a special left $\GInj\La$-approximation of $X$.
\end{enumerate}
\end{prop}
Note that the equivalence $\bfK_\tac(\Inj\La)\xto{\sim} \oGInj\La$
restricts to an equivalence $$\bfK_\tac(\inj\La)\stackrel{\sim}\lto
\oGInj\La\cap\omod\La.$$

\subsection{Torsion pairs for triangulated categories}

Let $\T$ be a triangulated category, for instance $\T=\bfK(\Inj\La)$.
A full subcategory $\U$ is {\em thick} if it is closed under shifts,
mapping cones, and direct factors. We denote by $\Thick(\U)$ the
smallest thick subcategory of $\T$ which contains $\U$, and
$\Thickc(\U)$ denotes the smallest thick subcategory of $\U$ which
contains $\U$ and is closed under small coproducts.

\begin{defn}
  A pair $(\U,\V)$ of thick subcategories of $\T$ forms a {\em torsion
    pair} for $\T$ if the following conditions are satisfied:
\begin{enumerate}
\item $\U =\{X\in\T\mid\Hom_\T(X,Y)=0\text{ for all }Y\in\V\}$;
\item $\V =\{Y\in\T\mid\Hom_\T(X,Y)=0\text{ for all }X\in\U\}$;
\item every object $X\in\T$ fits into an exact triangle $X'\to X\to
X''\to $ with $X'\in\U$ and $X''\in\V$.
\end{enumerate}
\end{defn}

\subsection{A torsion pair for $\bfK(\Inj\La)$}\label{se:tor}

The totally acyclic complexes of injective $\La$-modules form the
torsion free part of a torsion pair. This follows from the fact that a
complex $Y$ of injective $\La$-modules is totally acyclic if and only
if $$\Hom_{\bfK(\Inj\La)}(\bfi
X,Y[n])\cong\Hom_{\bfK(\Mod\La)}(X,Y[n])\cong H^n\Hom_\La(X,Y)=0$$ for
all $n\in\bbZ$ and all $X\in\Mod\La$ which are projective or
injective.

\begin{prop}[{\cite[\S 7]{K1}}]
We have a torsion pair
$$(\U,\V)=\big(\Thickc(\bfi(\proj\La\cup\inj\La)),\bfK_\tac(\Inj\La)\big)$$
for $\bfK(\Inj\La)$ which has the following properties.
\begin{enumerate}
\item The pair $(\U,\V)$ induces a cotorsion pair
$$(\bfs\U,\bfs\V)=({^\perp(\GInj\La)},\GInj\La)$$ for $\Mod\La$ with
  $${^\perp(\GInj\La)}\cap\GInj\La=\Inj\La.$$
\item Let $X\in\Mod\La$ and choose an exact triangle $X'\xto{\p} \bfi
X\xto{\psi} X''\to$ in $\bfK(\Inj\La)$ with $X'\in\U$ and
$X''\in\V$. Then $\bfs\p$ represents a special right
${^\perp(\GInj\La)}$-approximation and $\bfs\psi$ represents a special
left $\GInj\La$-approximation of $X$.
\end{enumerate}
\end{prop}

Let us mention the following dual result (which will not be used in
this work). There exists a torsion pair
$$\big(\bfK_\tac(\Proj\La),\Thickp(\bfp(\proj\La\cup\inj\La))\big)$$ for
$\bfK(\Proj\La)$ which induces via the canonical functor
$\bfK(\Proj\La)\to\uMod\La$ the cotorsion pair
$(\GProj\La,(\GProj\La)^\perp)$ for $\Mod\La$.

\begin{rem}
If a finitely generated $\La$-module $X$ admits a special left
$\GInj\La$-approximation in $\mod\La$, then $X$ admits a special right
$^\perp(\GInj\La)$-approximation in $\mod\La$. This follows from the
fact that both approximations are connected via an exact triangle
$X'\xto{} \bfi X\xto{} X''\to$ which can be chosen to lie in
$\bfK(\inj\La)$ if $\bfs X''$ is finitely generated.
\end{rem}

\subsection{The finitely generated orthogonal complement of 
$\GInj\La$}\label{se:fgcomp}

\begin{thm}
We have 
$$(\GProj\La)^\perp\cap\mod\La=\Thick(\proj\La\cup\inj\La)=
{^\perp(\GInj\La)}\cap\mod\La.$$
\end{thm}
\begin{proof}
Fix a finitely generated $\La$-module $X$. We use the cotorsion pair
$(\U,\V)$ from Proposition~\ref{se:tor}. Observe that $X$ belongs
to ${^\perp(\GInj\La)}$ if and only if $\bfi X$ belongs to $\U$. This
follows from part (2) of Proposition~\ref{se:tor} if we consider the
exact triangle $X'\xto{} \bfi X\xto{} X''\to$ with $X'\in\U$ and
$X''\in\V$. We have $X\in{^\perp(\GInj\La)}$ if and only if $\bfs
X''\in{^\perp(\GInj\La)}\cap\GInj\La=\Inj\La$.  The complex $X''$ is
totally acyclic and therefore Proposition~\ref{se:tac} implies that
$\bfs X''\in\Inj\La$ if and only if $X''=0$ if and only if $\bfi
X\in\U$. Next we use the identification
$$\bfD^b(\mod\La)=\bfK(\Inj\La)^c$$ via the functor $\bfi$ sending a
complex to its semi-injective resolution; see (\ref{se:kinj}).  This
implies
$$\U\cap\bfK(\Inj\La)^c=\Thick(\bfi(\proj\La\cup\inj\La))$$ and we see
that $X\in {^\perp(\GInj\La)}$ if and only if $\bfi
X\in\Thick(\bfi(\proj\La\cup\inj\La))$. Applying the correspondence
between thick subcategories of $\mod\La$ and $\bfD^b(\mod\La)$ from
the appendix, we conclude that $X\in {^\perp(\GInj\La)}$ if and only
if $X\in\Thick(\proj\La\cup\inj\La)$.

To complete the proof, observe that we have 
$$(\GProj\La)^\perp\cap\mod\La= {^\perp(\GInj\La)}\cap\mod\La$$ by
\cite[Lemma~8.6]{B}.
\end{proof}

\section{Virtually Gorenstein algebras}

\subsection{A characterization via thick subcategories}\label{se:char}

We apply the results from the preceding section and provide a
characterization of virtually Gorenstein algebras in terms of finitely
generated modules.
\begin{thm}
For an Artin algebra $\La$ the following are equivalent.
\begin{enumerate}
\item The algebra $\La$ is virtually Gorenstein.
\item The subcategory $\Thick(\proj\La\cup\inj\La)$ of $\mod\La$ is
contravariantly finite.
\item The subcategory $\Thick(\proj\La\cup\inj\La)$ of $\mod\La$ is
covariantly finite.
\end{enumerate}
\end{thm} 
\begin{proof}
  (1) $\Rightarrow$ (2): We adapt the proof of Theorem~8.2 in
  \cite{B}, using Theorem~\ref{se:fgcomp} as a new ingredient. In
  fact, the material which we need has been collected in
  (\ref{se:tac}) and (\ref{se:tor}). The proof is based on the
  following elementary fact. Given an adjoint pair of functors, the
  left adjoint preserves compactness provided the right adjoint
  preserves small coproducts.  Note that $(\uMod\La)^c=\umod\La$ and
  $(\oMod\La)^c=\omod\La$; see \cite[Lemma~6.3]{B}. We apply this fact
  to the inclusion $\uGProj\La\to\uMod\La$. This functor has a right
  adjoint $G\colon\uMod\La\to\uGProj\La$ such that the adjunction
  morphism $GX\to X$ represents a special right
  $\GProj\La$-approximation for every $\La$-module $X$; see
  (\ref{se:tac}).  The equality $(\GProj\La)^\perp={^\perp(\GInj\La)}$
  implies that $(\GProj\La)^\perp$ is closed under small coproducts.
  Thus a small coproduct of special right $\GProj\La$-approximation is
  again a special right $\GProj\La$-approximation. It follows that $G$
  preserves coproducts and therefore every compact object in
  $\uGProj\La$ belongs to $\umod\La$. The equivalence
  $\Proj\La\xto{\sim}\Inj\La$ induces an equivalence
$$\uMod\La\supseteq\uGProj\La\xleftarrow{\sim}
\bfK_\tac(\Proj\La)\stackrel{\sim}\to\bfK_\tac(\Inj\La)
\xto{\sim}\oGInj\La\subseteq\oMod\La$$ which sends compact objects to
compact objects and objects in $\umod\La$ to objects in
$\omod\La$. Thus every compact object in $\oGInj\La$ belongs to
$\omod\La$. On the other hand, the left adjoint
$F\colon\oMod\La\to\oGInj\La$ of the inclusion preserves compactness,
since the inclusion preserves small coproducts. Given $X\in\Mod\La$,
the adjunction morphism $X\to FX$ represents a special left
$\GInj\La$-approximation of $X$. We conclude that this approximation
can be chosen in $\mod\La$ if $X$ is finitely generated. The
Remark~\ref{se:tor} shows that in this case $X$ admits a special right
$^\perp(\GInj\La)$-approximation in $\mod\La$.  Thus
$$\Thick(\proj\La\cup\inj\La)={^\perp(\GInj\La)}\cap\mod\La$$ is a
contravariantly finite subcategory of $\mod\La$, thanks to
Theorem~\ref{se:fgcomp}.

(2) $\Leftrightarrow$ (3): In \cite[Corollary~2.6]{KS}, it is shown
    that every resolving and contravariantly finite subcategory of
    $\mod\La$ is covariantly finite. Dually, every coresolving and
    covariantly finite subcategory of $\mod\La$ is contravariantly
    finite.

    (2) \& (3) $\Rightarrow$ (1): Let
    $\D=\Thick(\proj\La\cup\inj\La)$.  We obtain for $\mod\La$ a
    cotorsion pair $(\C,\D)$ because $\D$ is a covariantly finite and
    coresolving subcategory. This is a well-known fact; see for
    instance \cite{AR} or \cite[Lemma~2.1]{KS}.  Analogously, we
    obtain a cotorsion pair $(\D,\E)$ for $\mod\La$ since $\D$ is a
    contravariantly finite and resolving subcategory.  Given any
    subcategory $\X$ of $\mod\La$, we denote by $\li\X$ the full
    subcategory of $\Mod\La$ consisting of all filtered colimits of
    modules in $\X$. In \cite[Theorem~2.4]{KS}, it is shown that there
    are cotorsion pairs $(\li\C,\li\D)$ and $(\li\D,\li\E)$ for
    $\Mod\La$. We claim that
$$(\li\C,\li\D)=(\GProj\La,(\GProj\La)^\perp) \quad\text{and} \quad
(\li\D,\li\E)=({^\perp(\GInj\La)},\GInj\La).$$ To prove this claim,
first observe that $\li\D$ is resolving and coresolving. Therefore
$\li\C\subseteq\GProj\La$ and $\li\E\subseteq\GInj\La$ by
Lemma~\ref{se:resolv}.  Now fix a special right $\D$-approximation
$S_\D\to S$ of $S=\La/{\rad\La}$ and observe that $S_\D$ belongs to
${^\perp(\GInj\La)}$ since $\D\subseteq {^\perp(\GInj\La)}$ by
Theorem~\ref{se:fgcomp}. This implies
$\li\D\subseteq{^\perp(\GInj\La)}$ because each object in $\li\D$ is
obtained from $S_\D$ by taking small coproducts of copies of $S_\D$,
forming finitely many extensions, and taking direct factors; see
\cite[Theorem~2.4]{KS}. Analogously, $\li\D\subseteq
(\GProj\La)^\perp$ because each object in $\li\D$ is obtained from a
special left $\D$-approximation $S\to S^\D$ by taking small products
of copies of $S^\D$, finitely many extensions, and direct
factors. Thus our claim follows and therefore $\La$ is virtually
Gorenstein.
\end{proof}

Let $\La$ be a virtually Gorenstein algebra. Then we obtain from the
preceding proof a description of the subcategory
$(\GProj\La)^\perp={^\perp(\GInj\La)}$. To formulate this, we use the
notation $\Thickp$ and $\Thickc$ to denote thick subcategories of
$\Mod\La$ which are closed under small products and coproducts,
respectively.

\begin{cor}
Let $\La$ be a virtually Gorenstein algebra. Then we have
$$(\GProj\La)^\perp=\Thickp(\Proj\La\cup\Inj\La)=\Thickc(\Proj\La\cup\Inj\La)=
{^\perp(\GInj\La)}.$$
\end{cor}

\subsection{A characterization via filtered colimits}\label{se:char2}

We provide another characterization of virtually Gorenstein algebras
in terms of filtered colimits of finitely generated modules. It is
convenient to define the subcategories
$$\Gproj\La:=\GProj\La\cap\mod\La\quad\text{and}\quad
\Ginj\La:=\GInj\La\cap\mod\La.$$

\begin{thm} 
For an Artin algebra $\La$ the following are equivalent. 
\begin{enumerate}
\item $\Lambda$ is virtually Gorenstein.  
\item Any Gorenstein projective module is a filtered colimit of
finitely generated Gorenstein projective modules.
\item Any Gorenstein injective module is a filtered colimit of
finitely generated Gorenstein injective modules.
\end{enumerate}
If $\Lambda$ is virtually Gorenstein, then $\Gproj\La$ and $\Ginj\La$
are both covariantly and contravariantly  finite subcategories of $\mod\La$.
\end{thm}
\begin{proof}
(1) $\Rightarrow$ (2) $\&$ (3): See the proof of
    Theorem~\ref{se:char}.

(2) $\Leftrightarrow$ (3): Let $D$ denote the duality between right
and left $\La$-modules.  The adjoint pair of functors
$-\otimes_{\La}D(\Lambda)$ and $\Hom_{\La}(D(\Lambda),-)$ induces an
equivalence between $\Proj\La$ and $\Inj\La$, which restricts to an
equivalence between $\proj\La$ and $\inj\La$. Therefore
the adjoint pair induces an equivalence between $\GProj\La$ and
$\GInj\La$, and between $\Gproj\La$ and $\Ginj\La$, respectively, see
\cite[Proposition~3.4]{B}. Now use the fact that both functors
preserves filtered colimits.

(3) $\Rightarrow$ (1): A standard argument shows that $\Ginj\La$ is a
covariantly finite subcategory of $\mod\La$. To see this, fix $X$ in
$\mod\La$ and let $X\to Y_\a$ be a representative family of maps into
objects from $\Ginj\La$. Then $\prod_\a Y_\a$ is by our assumption a
filtered colimit of objects in $\Ginj\La$, and therefore the map $X\to
\prod_\a Y_\a$ factors through a map $ X\to Y$ with $Y$ in
$\Ginj\La$. By our construction, the map $X\to Y$ is a left
$\Ginj\La$-approximation of $X$, since every map $X\to Y_\a$ factors
through $X\to Y$. Thus $\Ginj\La$ is covariantly finite.

Let $\D=\Ginj\La$. As in the proof of Theorem~\ref{se:char}, we obtain
cotorsion pairs $(\C,\D)$ and $(\li\C,\li\D)$ for $\mod\La$ and
$\Mod\La$, respectively. We claim that
$\C={^\perp(\GInj\La)}\cap\mod\La$. We have $\C\subseteq
{^\perp(\GInj\La)}\cap\mod\La$, because $\GInj\La\subseteq\li\D$ by
our assumption. To show the other inclusion, let $X$ be in
${^\perp(\GInj\La)}\cap\mod\La$ and choose an exact sequence $0\to
D\to C\to X\to 0$ such that $C\to X$ is a special right
$\C$-approximation. This sequence splits since $D$ is Gorenstein
injective, and therefore $X$ belongs to $\C$. It remains to recall
from Theorem~\ref{se:fgcomp} that
$${^\perp(\GInj\La)}\cap\mod\La=\Thick(\proj\La\cup\inj\La).$$ Thus
$\Thick(\proj\La\cup\inj\La)$ is contravariantly finite, and therefore the
characterization from Theorem~\ref{se:char} implies that $\La$ is
virtually Gorenstein.

To conclude this proof, assume that $\La$ is virtually
Gorenstein. Thus $\La^\op$ is virtually Gorenstein as well. We have
seen that $\Ginj\La$ is a covariantly finite subcategory of
$\mod\La$. Recall that any covariantly finite and coresolving
subcategory of $\mod\La$ is contravariantly finite, by
\cite[Corollary~2.6]{KS}. Thus $\Ginj\La$ is contravariantly
finite. The duality between right and left $\La$-modules identifies
$\Ginj\La^\op$ with $\Gproj\La$. We conclude that $\Gproj\La$ is both
covariantly and contravariantly finite.
\end{proof}

\subsection{An example}
We provide an example of an Artin algebra which is not virtually
Gorenstein.  This is based on work of Yoshino \cite{Y}. In fact, he
constructs a class of algebras $\La$ such that the subcategory
$\Gproj\La$ of $\mod\La$ is not contravariantly finite. Then we use
the fact from Theorem~\ref{se:char2} that $\Gproj\La$ is
contravariantly finite whenever $\La$ is virtually Gorenstein.

\begin{prop}
Let $K$ be a field. Then the $6$-dimensional $K$-algebra $$\Lambda =
 K[x,y,z]/\langle x^{2},yz,y^{2}-xz,z^{2}-yx \rangle$$
is not virtually Gorenstein.
\end{prop}
\begin{proof}
Consider the $K$-algebra
$$\Gamma = K[x,y,z]/ \langle xz-y^{2}, yx-z^{2}, zy-x^{2}\rangle.$$
This is a one-dimensional Cohen-Macaulay non-Gorenstein homogeneous
ring. We have $\Lambda = \Gamma/x^{2}\Gamma$ and this algebra has
radical cubed zero and is non-Gorenstein. It is shown in \cite{Y} (in
a more general setting) that the trivial $\Lambda$-module $K$ has no
right $\Gproj\La$-approximation, hence $\Gproj\La$ fails to be
contravariantly finite in $\mod\La$. The proof uses the graded
structure of $\La$ and its Hilbert series. Consequently, by
Theorem~\ref{se:char2}, $\La$ is an example of a finite-dimensional
$K$-algebra which is not virtually Gorenstein.
\end{proof}

\begin{appendix}
\section{Thick subcategories}

Let $\A$ be an exact category and suppose that $\A$ is idempotent
complete. A full subcategory $\C$ of $\A$ is called {\em thick} if it
is closed under direct factors and has the following two out of three
property: for every exact sequence $0\to X\to Y\to Z\to 0$ in $\A$
with two terms in $\C$, the third term belongs to $\C$ as well.

Now consider the bounded derived category $\bfD^b(\A)$ of $\A$ and
identify $\A$ with the full subcategory of $\bfD^b(\A)$ consisting of
all complexes concentrated in degree zero. Recall that a full
subcategory of a triangulated category is {\em thick} if it is closed
under shifts, mapping cones, and direct factors. 

We discuss the relation between thick subcategories of exact and
triangulated categories. To this end, let us call a thick subcategory
$\C$ of $\A$ {\em cofinal} if every admissable epimorphism $Y\to Z$
into an object $Z\in\C$ admits an admissable epimorphism $Y'\to Z$
with $Y'\in\C$ and factoring through $Y\to Z$.

\begin{prop}
Let $\C$ be a full subcategory of $\A$.  If $\C$ is a cofinal thick
 subcategory of $\A$, then $\C$ is of the form $\D\cap\A$ for some
 thick subcategory $\D$ of $\bfD^b(\A)$.  Conversely, if
 $\C=\D\cap\A$ for some thick subcategory $\D$ of $\bfD^b(\A)$, then
 $\C$ is a thick subcategory of $\A$.
\end{prop}
\begin{proof}
Suppose that $\C$ is a cofinal thick subcategory.  This assumption on
$\C$ implies that the inclusion $\C\to\A$ induces a fully faithful and
exact functor $\bfD^b(\C)\to\bfD^b(\A)$; see for instance
\cite[Proposition~III.2.4.1]{V}. Note also that idempotents in
$\bfD^b(\C)$ split since $\C$ has this property; see
\cite[Theorem~2.8]{BS}.  Thus the full subcategory $\D$ of
$\bfD^b(\A)$ consisting of complexes quasi-isomorphic to a complex of
objects in $\C$ is a thick subcategory. We claim that
$\C=\D\cap\A$. Clearly, $\C\subseteq\D\cap\A$. Thus we fix
$X\in\D\cap\A$. Then $X$ is in $\bfD^b(\A)$ isomorphic to a bounded
complex $C$ with differential $\d$ such that $C^n\in\C$ for all $n$
and $C$ is acyclic in all degrees $n\neq 0$. Now we use that $\C$ is
thick. Thus $\Coker\d^{-2}$ and $\Ker\d^0$ belong to $\C$, and we have
an admissable monomorphism $\Coker\d^{-2}\to\Ker\d^0$ such that the
cokernel is isomorphic to $X$. We conclude that $X$ belongs to $\C$
and therefore $\C=\D\cap\A$.

The converse is clear since each exact sequence $0\to X\to Y\to Z\to
0$ in $\A$ gives rise to an exact triangle $X\to Y\to Z\to$ in
$\bfD^b(\A)$.
\end{proof}

\begin{exm}
Let $\A$ be an exact category having enough projective objects.  Then
every thick subcategory $\C$ containing all projective objects is cofinal.
\end{exm}

\subsection*{Acknowledgement}
The example of an algebra which is not virtually Gorenstein arose from
discussions during a workshop on thick subcategories at Oberwolfach in
February 2006. We are grateful to Osamu Iyama, Srikanth Iyengar and
Apostolos Thoma.

\end{appendix}

\end{document}